\documentclass[12pt]{amsart}
\usepackage{amssymb}
\usepackage{graphicx}

\newtheorem{theorem}{Theorem}[section]

\newtheorem{proposition}[theorem]{Proposition}
\newtheorem{corollary}[theorem]{Corollary}
\newtheorem{lemma}[theorem]{Lemma}

\def \proof {\noindent {\bf Proof.}\ \ }
\def \qed {{\mbox{}\nolinebreak\hfill\rule{2mm}{2mm}\par\medbreak} }

\def \ll {\langle}
\def \rr {\rangle}
\def \R {{\bf R}}
\def \a {\alpha}

\def \e {\varepsilon}
\def \d {\delta}

\def \l {\lambda}
\def \L {\Lambda}
\def \la {\lambda}

\def \R {\mathbb{R}}
\def \< {\langle}
\def \> {\rangle}

\def \inte {{\rm int}}

\def \Z {\mathbb{Z}}
\def \C {\mathbb{C}}
\def \N {\mathbb{N}}

\begin{document}
\title {Riesz basis of exponentials for a union of cubes in $\R^{d}$}
\author {Jordi Marzo}
\address{Departament de matem\`atica aplicada i an\`alisi, Universitat de Barcelona,
                               Gran via 585, 08071 Barcelona, Spain}
\email{jmarzo@mat.ub.es}

\thanks{The author is supported by the DGICYT grant: BES-2003-2618, MTM2005-08984-C02-02
.}



\date{December 15, 2005}

\begin{abstract}
  We extend to several
  dimensions the result of K.~Seip and Y.I.~Lyubarskii
  that proves the existence of Riesz basis of exponentials
  for a finite union of intervals with equals lengths.
\end{abstract}

\maketitle

\section{Introduction}                                      \label{SecIntro}

Let $E \subset \R^{d}$ be a bounded set with no empty interior. We
consider the Paley-Wiener space of $L^{2}$ functions band-limited to
$E$
$$PW_{E}=\{ F\in L^{2}(\R^{d}):\sup \mathcal{F}F\subset E\},$$
endowed with the $L^{2}(\R^{d})$ norm. Here we denote as
$$\mathcal{F}F(x)=\frac{1}{(2\pi)^{d/2}}\int_{\R^{d}}F(t)e^{-ix\cdot t}dt,$$
the Fourier transform of $F.$

By the Paley-Wiener theorem $PW_{E}$ is a reproducing kernel Hilbert
space of entire functions with kernel
$$K(z,w)=\frac{1}{(2\pi)^{d/2}}\int_{E}e^{it\cdot (z-\overline{w})}dt.$$
Then, for any $F\in PW_{E},$
$$F(w)=\frac{1}{(2\pi)^{d/2}}\int_{\R^{d}}F(x)K(w,x)dx
,\;\;\; w\in \C^{d}.$$

We say that a sequence $\Lambda=\{ \lambda_{n}\} \subset\R^{d}$ is
\emph{sampling} for $PW_{E}$
    if there exist constants $A,B>0$ such that
\begin{equation}                                                            \label{sampling}
    A \, \int_{\R^{d}}|F(x)|^{2}dx \leq \sum_{n}|F(\l_{n})|^{2}\le B \,
    \int_{\R^{d}}|F(x)|^{2}dx,
\end{equation}
    for all $F\in PW_{E}.$

We say that a sequence $\Lambda=\{ \lambda_{n}\}
\subset\mathbb{R}^{d}$ is \emph{interpolating} for $PW_{E}$ if for
    each
    $\{a_{n}\}\in\ell^{2}$ there exists a function
    $F\in PW_{E}$ such that
\begin{equation}                                                            \label{interpolation}
    F(\la_{n})=a_{n}, \quad n\in \mathbb{Z}.
\end{equation}

If $\Lambda$ is both sampling and interpolating we say that it is a
\emph{complete interpolating sequence.}

Sampling and interpolation sequences in Paley-Wiener spaces provide
a mathematical model of stable recovery and data transmission in
signal theory regardless of the computation method. For the meaning
of this concepts we refer to \cite{L}. In 1967 H. J. Landau found
necessary conditions for the existence of such a sequences in terms
of some densities, as a consequence a complete interpolating
sequence must have density equal to the Landau-Nyquist rate, see
Theorem \ref{lan}.

The aim of this paper is to find a complete interpolating sequence
for $PW_{E},$ where $E$ is a finite union of mutually disjoint
half-open cubes
\begin{equation}                                                            \label{conj}
E=\bigcup_{j=1}^{p} Q_{j}\subset\mathbb{R}^{d},\;\; \mbox{where}
\;\; Q_{j}=\prod_{s=1}^{d}[\alpha_{j}^{s},\alpha_{j}^{s}+\beta),\;\;
\beta>0.
\end{equation}
Our result allows us to construct, for any bounded $E$ with no empty
interior, sampling or interpolating sequences with density
arbitrarily close to the critical one, see Corollary \ref{coro}.
This implies that Landau's results \cite{LAN} about densities can't
be improved.

Moreover, we find a recovery process in $PW_{E}$ whose
stability is assured by the sampling property, see Theorem
\ref{teo:2}. This recovery process can be considered as a shifted version of the
Shannon-Whittacker-Kotelnikov theorem.

The method used is a generalization to $\R^{d}$ with
$d>1$ to that of \cite{LS} for $\R.$ It is worth noting that the
reconstruction formula (\ref{eq-fin}), which comes from A.
Kohlenberg's work \cite{K}, has been used by A. Faridani to produce
sampling schemes for computerized tomography \cite{F}.
Our result
provides the stability that their result lacks. In fact, Faridani
proved a more general version of the
reconstruction formula (\ref{eq-fin}) that holds for smooth
functions with spectrum in any compact set of $\R^{d},$
and there is a version even for the group context \cite{BF}.
But in such
a general context there is no hope of stability.

Finally, we observe that we can define for $1< q<\infty$ the
Paley-Wiener space $PW_{E}^{q}$ of entire functions given by
$$F(z)=\frac{1}{(2\pi)^{d/2}}\int_{E}f(t)e^{it\cdot z}dt, \;\;z\in \C^{d}$$
with $f\in L^{q^{'}}(E),$ for $q^{'}=q/(q-1).$ With the natural modifications
the main result (Theorem \ref{teo:2})
holds also in $PW_{E}^{q}$ for $1<q<\infty$ because the
projection of $L^{q}(\R^{d})$ into $PW^{q}(E)$ is bounded and
we have also a Plancherel-P\'olya type inequality \cite[3.2.2.]{NIK}.

\section{Preliminaries}                                                \label{SecPrelim}

In this section we review some results and definitions that will be used later on.

Let $v=\{ v_{n} \}$ be a sequence in a separable Hilbert space $H.$
We say that $v$ is a \emph{frame} for $H$ if there exist constants $A,B>0$
such that
$$A|| u ||^{2}\le \sum_{n}|\ll u,v_{n} \rr|^{2}\le B|| u ||^{2},\;\; u\in H.$$
We say that $v=\{ v_{n} \}$ is a \emph{Riesz sequence} or, equivalently, that
it is a Riesz basis in its closed linear span, if there exist
constants $A,B>0$ such that for every finite scalar sequence $\{
c_{n} \}$ one has
$$A\sum|c_{n}|^{2}\le \left\| \sum c_{n}v_{n} \right\|^{2} \le B\sum|c_{n}|^{2}.$$

Now we can translate our problem concerning sequences into the basis language.
A sequence $\L=\{\l_{n} \} \subset \R^{d}$ is interpolating
for $PW_{E}$ if and only if the sequence $\{ K(\cdot
,\l_{n})/\sqrt{K(\l_{n},\l_{n})} \}$ is a Riesz sequence for
$PW_{E},$ and this is also equivalent to the assumption that
$$\mathcal{E}(\L)=\{ e^{i \l_{m} \cdot x}: \l_{m} \in \L \}$$
is a Riesz sequence for $L^{2}(E).$ Similarly, $\L=\{\l_{n}
\} \subset \R^{d}$ is sampling if and
only if $\{ K(\cdot ,\l_{n})/\sqrt{K(\l_{n},\l_{n})}
\}$ is a frame for $PW_{E},$ and
if and
only if
$$\mathcal{E}(\L)=\{ e^{i \l_{m} \cdot x}: \l_{m} \in \L \}$$
is a frame for $L^{2}(E).$
Then $\L$ is both
sampling and interpolating if and only if $\mathcal{E}(\L)$ is a
Riesz basis for $L^{2}(E)$, i.e. $\mathcal{E}(\L)$ is the image of
an orthonormal basis by a bounded invertible operator, see
\cite{LS}.

We say that a sequence $\Lambda=\{ \l_{n}
    \} \subset\R^{d}$ is \emph{uniformly discrete} if
    $$\inf_{j\neq i} |\la_{j}-\la_{k}|>0.$$
This infimum is called the separation constant of $\L.$

For $h>0$ and
$x=(x_{1},\dots ,x_{d})$ denote by $Q_{h}(x)$ the cube:
$$Q_{h}(x)=\prod_{s=1}^{d}\left[x_{s}-h/2,x_{s}+h/2\right).$$
Then $\{ Q_{h}(hm) \}_{m\in \Z^{d}}$ is a disjoint cover of $\R^{d}$ for any $h>0.$

Let $\L=\{ \l_{n} \}\subset \R^{d}$ be a sequence.
Following Beurling we define the upper uniform density $D^{+}(\L)$ and lower uniform density $D^{-}(\L)$
of $\L$ as follows:
\begin{equation*}
    D^{+}(\Lambda)=\limsup_{h\to \infty}\frac{n^{+}(h)}{h^{n}}\;\; \mbox{and} \;\;
            D^{-}(\Lambda)=\liminf_{h\to
            \infty}\frac{n^{-}(h)}{h^{n}},
\end{equation*}
where
\begin{equation*}
    n^{+}(h)=\sup_{x\in \R^{d}} \#(\Lambda\cap Q_{h}(x))\;\; \mbox{and}\;\;
    n^{-}(h)=\inf_{x\in \R^{d}} \#(\Lambda\cap Q_{h}(x)).
\end{equation*}
If this densities coincide we say
that $\L$ has uniform density and we write $D(\L)=D^{+}(\L)=D^{-}(\L).$

The following result, which provides necessary conditions for
sampling and interpolation, is due to H. J. Landau \cite{LAN}.
\begin{theorem}                                                                 \label{lan}
    If
    $\L$ is an interpolating sequence for $PW_{E}$ then $\L$ is uniformly
    discrete and
    $$D^{+}(\L)\le \frac{|E|}{(2\pi)^{d}}.$$
    If $\L$ is a sampling sequence for $PW_{E}$ there exists
    $\L^{'}\subset \L$ uniformly discrete and sampling subsequence such
    that
    $$D^{-}(\L^{'})\ge \frac{|E|}{(2\pi)^{d}}.$$
\end{theorem}

Then Riesz basis can occur only when $D(\L)=|E|/(2\pi)^{d}.$ This particular value
is called
the Landau-Nyquist rate.

It is quite easy to show that no characterization can be obtained
for sampling and interpolation in several dimensions, nor in the
one-dimensional multi-band case, using only these densities.

Unlike the case of a single interval, where there is a complete
characterization for the complete interpolating sequences \cite{Pa},
little is known for the same problem with several intervals or in
several dimensions. Let us comment some of the known results. The
reasonable question in this more general setting
 is to find some complete interpolation sequence
rather than to obtain a complete characterization. For a
finite union of commensurable intervals, for two general intervals or for
particular configurations with more than two intervals, it is known that such
sequences exist \cite{LS}. In $\R^{2}$ there are also complete
interpolating sequences for $PW_{E}$ if $E$ is a convex polygon
symmetric with respect to the origin \cite{LR}. In $\R^{d}$ for $d>2$ the only examples,
as far as we know, concern the spaces $PW_{E}$ with parallelogram domains $E.$
For more on this
problem see \cite{Se}.

Our result provides the first example of non-convex and
non-symmetric set with Riesz basis of exponentials in several
dimensions. This result contrast with the more rigid analogue problem of
orthonormal basis of exponentials for $L^{2}(E),$ where it is
known that symmetry is necessary \cite{Kolo}.

Our result
allows us to construct, for any $E\subset \R^{d}$ bounded with no
empty interior, sampling or interpolating sequences for $PW_{E}$
with density arbitrarily close to the Landau-Nyquist rate, see
Corollary \ref{coro}. Unfortunately this method doesn't seem to be
useful to find complete interpolating sequences for $PW_{E}$ with
$E$ arbitrary.

Next we will show that the separation property of a sequence is
equivalent to the right inequality in (\ref{sampling}). This is a
well known result but we include its proof for the sake of
completeness \cite[pp 82-83]{Y}.
\begin{lemma}                                                          \label{lem:1}
    Let $E\subset \R^{d}$ be a bounded set with no empty interior and
    $\L=\{ \l_{n} \} \subset \R^{d}$ be a sequence.
   $\L$
    is a finite union of uniformly discrete sequences if and only if
     there exists a constant $B>0$ such that
\begin{equation}                                                       \label{ineq:1}
    \sum_{n}\left| F(\l_{n})
    \right|^{2}\le
    B \int_{\R^{d}} |F(x)|^{2} dx,
\end{equation}
    for all $F\in PW_{E}.$
    Moreover, the constant $B>0$ depends only
    on the diameter of $E$ and on the separation
    constant of $\L$.
\end{lemma}

\vskip .2in

\proof Suppose, to reach a contradiction, that we have (\ref{ineq:1}) but
$\L$ is not a finite union of uniformly discrete sequences. We can suppose that $0\in
\inte \,E$ and then there exists an $\e>0$ such that $Q_{\e}(0)\subset
E$ and
$$\sup_{m\in \Z^{d}}\# (\L\cap Q_{h}(hm))=\infty,$$
where $h=\frac{\pi}{2\e}.$
Choose a sequence $\{m_{j}\}_{j\in\Z}\subset \Z^{d}$ such that
$$\# (\L\cap Q_{h}(hm_{j}))\to \infty,\;\; j\to \infty.$$
The functions $f_{j}(x)=e^{-ih m_{j}\cdot x}\chi_{Q_{\e}(0)}(x)$
have support on $E$ and
\begin{eqnarray*}
 F_{j}(z) & = & \mathcal{F}^{-1}f_{j}(z)=\frac{1}{(2\pi)^{d/2}}\int_{\R^{d}}f_{j}(x)e^{iz\cdot x} dx
 \\
 & = & \frac{1}{(2\pi)^{d/2}}\prod_{s=1}^{d}
\frac{\sin \e(z^{s}-h m_{j}^{s})}{z^{s}-h m_{j}^{s}},
\end{eqnarray*}
are in $PW_{E}.$
Now
$$\inf_{j\in \Z}\inf_{z\in Q_{h}(hm_{j})}|F_{j}(z)|>0,$$
because
\begin{eqnarray*}
    \inf_{z\in Q_{h}(hm_{j})}|F_{j}(z)| & = & \frac{1}{(2\pi)^{d/2}}\inf_{z\in Q_{h}(0)}\left|
    \prod_{s=1}^{d}\frac{\sin \e z^{s}}{z^{s}} \right| \\
    & \ge &
    \left( \frac{\e}{\sqrt{2\pi}}\inf_{|t|<1/2}\left|\frac{\sin \pi t}{\pi t}\right|  \right)^{d},
\end{eqnarray*}
and therefore
$$||F_{j}|\L||_{\ell^{2}}\to \infty,\;\; j\to \infty,$$
but $|| F_{j} ||=|Q_{\e}(0)|^{1/2}.$

To prove the converse, we can suppose that $\L$ is uniformly discrete,
i.e. there exists $\tau>0$ such that
$$|\l_{j}-\l_{k}|\ge \tau>0,\; \; j\neq k.$$
Take an interval that contains $E$
$$E\subset \prod_{s=1}^{d}[-a^{s},a^{s}].$$
Then $F\in PW_{E}$ is of exponential type $\nu=(\nu^{1},\dots ,\nu^{d}),$
where $0<
\nu^{s}\le a^{s},$ $s=1,\dots ,p.$
From the mean-value property one has, for all
$\rho_{s}>0,$
$$F(x)=\frac{1}{(2\pi)^{d}}\int_{[0,2\pi]^{d}}
F(x^{1}+\rho_{1}e^{i\theta_{1}},\dots ,x^{d}+\rho_{d}e^{i\theta_{d}})d\theta, \;\; x\in \R^{d}.$$
Integrating
on $[0,\delta]^{d},$ for some $\delta>0,$
$$F(x)\left( \frac{\delta^{2}}{2}\right)^{d}=\frac{1}{(2\pi)^{d}}
\int_{\prod_{s=1}^{d}B(x^{s},\d)}F(z)dm(z),$$
where $B(x^{s},\d)=\{z\in\C : |x^{s}-z|<\d\}.$
Then by Cauchy-Schwarz inequality
\begin{equation*}
    |F(x) |^{2}\le
    B
    \int_{[-\d,\d]^{d}}\left(
    \int_{\prod_{s=1}^{d}Q_{2\d}(x^{s})}
     |
    F(\xi+i\eta) |^{2}d\xi
    \right) d\eta,
\end{equation*}
where $B>0$ depends on $\delta$ and $d.$

For all $F$ entire function of exponential type $(\nu^{1},\dots \nu^{d})$
such that $F \in L^{2}(\R^{d}),$ we have
the following inequality
of Plancherel-P\'olya type
\begin{equation}                                                                        \label{plan}
\| F(\cdot + i\eta) \|_{L^{2}(\R^{d})}\le e^{\sum_{s=1}^{d}\nu_{s}|\eta^{s}|} \| F \|_{L^{2}(\R^{d})},
\end{equation}
see \cite[3.2.2.]{NIK}.

Finally, if we take $\delta=\frac{\tau}{2}$ the cubes $Q_{2\delta}(\l_{j})$ are pairwise disjoint.
This together with
(\ref{plan}), gives
\begin{eqnarray*}
    \sum_{j}\left| F(\l_{j})\right|^{2}& \le &
    B\int_{[-\delta,\delta]^{d}}\sum_{j}
    \int_{\prod_{s=1}^{d}Q_{2\d}(\l_{j})}
    |
    F(\xi+i\eta) |^{2}d\xi
    d\eta \\
    & \le &
    B\int_{[-\delta,\delta]^{d}}
    \left(
    \int_{\R^{d}}
    |
    F(\xi+i\eta) |^{2}d\xi
    \right)
    d\eta
    \\
    & \le &
    B\| F \|_{L^{2}(\R^{d})}^{2}
    \int_{[-\delta,\delta]^{d}}
    e^{\sum_{s=1}^{d}\nu_{s}|\eta^{s}|}
    d\eta \\
    & = & B\| F \|_{L^{2}(\R^{d})}^{2}.
\end{eqnarray*}
\qed


Following Kohlenberg, \cite{K}, we define the \emph{sampling procedure}
of order $p$ for a function $F\in PW_{E}$  such that
$\mathcal{F}F\in \mathcal{C}_{c}^{\infty}(E),$
as follows
\begin{equation}                                                                \label{samplingproc}
    G(t)
    =
    \sum_{l=1}^{p}\sum_{n\in \Z^{d}}
    F(W n+k_{l})
    S_{l}(t-Wn-k_{l}).
 \end{equation}
Here $W$ is a real non-singular $d\times d$ matrix,
$k_{l}\in \R^{d},$ and
$S_{l}\in PW_{E}.$
Since $\{ Wn+k_{l} \}_{n\in \Z^{d}}$ is uniformly discrete
(when the $k_{l}$ are different)
and $F$ tends to zero faster than any polynomial
this series converges.
We want to determine $W,$ $k_{l}$ and $S_{l}$ in order to get $G=F$ for all $F\in PW_{E}.$

First we compute the Fourier transform of this sampling procedure.
\begin{proposition}                                                                            \label{teo:1}
    Let $F,S_{l}\in PW_{E},$ $f=\mathcal{F}
    F\in \mathcal{C}^{\infty}_{c}(E),$
    $s_{l}=\mathcal{F} S_{l}$ and $G$ given by (\ref{samplingproc}).
    Then
\begin{equation}                                                                               \label{eq:1}
    \mathcal{F}G(\omega)= \frac{(2\pi)^{d/2}}{|\det W|}\sum_{l=1}^{p}
    s_{l}(\omega)
    \sum_{n \in\mathbb{Z}^{d}}
    f(\omega+2\pi (W^{-1})^{t}n)
    e^{2\pi i (W^{-1})^{t}m\cdot k_{l}}.
\end{equation}
\end{proposition}

\vskip .2in

\proof
    This result follows easily from Poisson's Formula and it can be
    found in \cite{F}.
    Consider the Fourier transform of G
\begin{equation*}
    \mathcal{F}G (\omega)=\sum_{l=1}^{p}
    s_{l}
    (\omega)
    \sum_{n \in\mathbb{Z}^{d}}
    F(Wn+k_{l})
    e^{-i\omega\cdot (Wn+k_{l})}.
\end{equation*}
    Poisson's formula gives
    $$(2\pi)^{d/2}|\det W|\sum_{n\in \Z^{d}}F_{k_{l}}(Wn)e^{-i\omega\cdot Wn}=
    \sum_{n\in \Z^{d}}\mathcal{F} F_{k_{l}}(\omega-2\pi (W^{-1})^{t}n)$$
    where $F_{k_{l}}=F(\cdot +k_{l}).$
\qed

\section{Construction and main result}                                               \label{SecMain}

Let $E$ be a finite union of half-open cubes
\begin{equation*}
E=\bigcup_{j=1}^{p} Q_{j}\subset\mathbb{R}^{d},\;\; \mbox{where}
\;\; Q_{j}=\prod_{s=1}^{d}[\alpha_{j}^{s},\alpha_{j}^{s}+\beta),\;\;
\beta>0.
\end{equation*}
We want to determine the functions and shifts $S_{l}(\cdot-k_{l}),$ $l=1,\dots ,p$
in order that the reconstruction formula
$$
F(t)=
    \sum_{l=1}^{p}\sum_{n\in \Z^{d}}
    F(W n+k_{l})
    S_{l}(t-Wn-k_{l})
$$
holds for all $F\in PW_{E}$ such that $f=\mathcal{F}F\in \mathcal{C}^{\infty}_{c}(E).$
We know, by Proposition \ref{teo:1}, that for such functions this is equivalent to
\begin{equation}                                                                                \label{marca}
f(\omega)= \frac{(2\pi)^{d/2}}{|\det W|}\sum_{l=1}^{p}
    s_{l}(\omega)
    \sum_{n \in\mathbb{Z}^{d}}
    f(\omega+2\pi (W^{-1})^{t}n)
    e^{2\pi i (W^{-1})^{t}n\cdot k_{l}}.
\end{equation}
From now on we take $W=\frac{2\pi}{\beta}Id$ and,
as $f=0$ off of $E$ we define $s_{l}(\omega)=0$ if $\omega \not \in E.$ To determine $s_{l}$ on $E$
we split $E$ in the parts which are overlapped by one or more of the sets
$$E+2\pi (W^{-1})^{t}n,\;\; n\in \Z^{d}.$$ Let us see how this partition works.

Consider two cubes $Q_{j}$ and $Q_{k}.$ There exists a
unique $n_{jk}\in \Z^{d}$ such that
$$\gamma_{jk}=\a_{k}-\beta n_{jk}\in Q_{j},$$
where $\a_{k}=(\a_{k}^{1},\dots , \a_{k}^{d})$ is the lower left
point of $Q_{k}.$ Observe that $n_{jj}=0.$ Fixed $j\in \{1,\dots
,p\}$ split $Q_{j}$ into the half-open intervals formed by
intersection of $Q_{j}$ with the $2^{d}$ quadrants around
$\gamma_{jk},$ and do the same for all of $k=1,\dots ,p.$ Thus, we
split each cube $Q_{j}$ in this way in at most $p^{d}$ intervals
(this can be easily proved by induction, projecting a face of the
cube into $\R^{d-1}$) which we denote $Q_{j}^{s}$ $s=1,\dots
,p^{d},$ the subscript indicates the cube in which the rectangle is
contained. From now on we assume that we have $p^{d}$ intervals
$Q_{j}^{s}$ on each cube $Q_{j}.$ The minor changes needed if that
is not the case are left to the reader. We have for each $j=1,\dots
,p$
$$\bigcup_{s=1}^{p^{d}} Q_{j}^{s}=Q_{j}, \;\;  Q_{j}^{s}\cap  Q_{j}^{\ell}=\emptyset.$$

\begin{figure}
\begin{center}
 \includegraphics[scale=0.7]{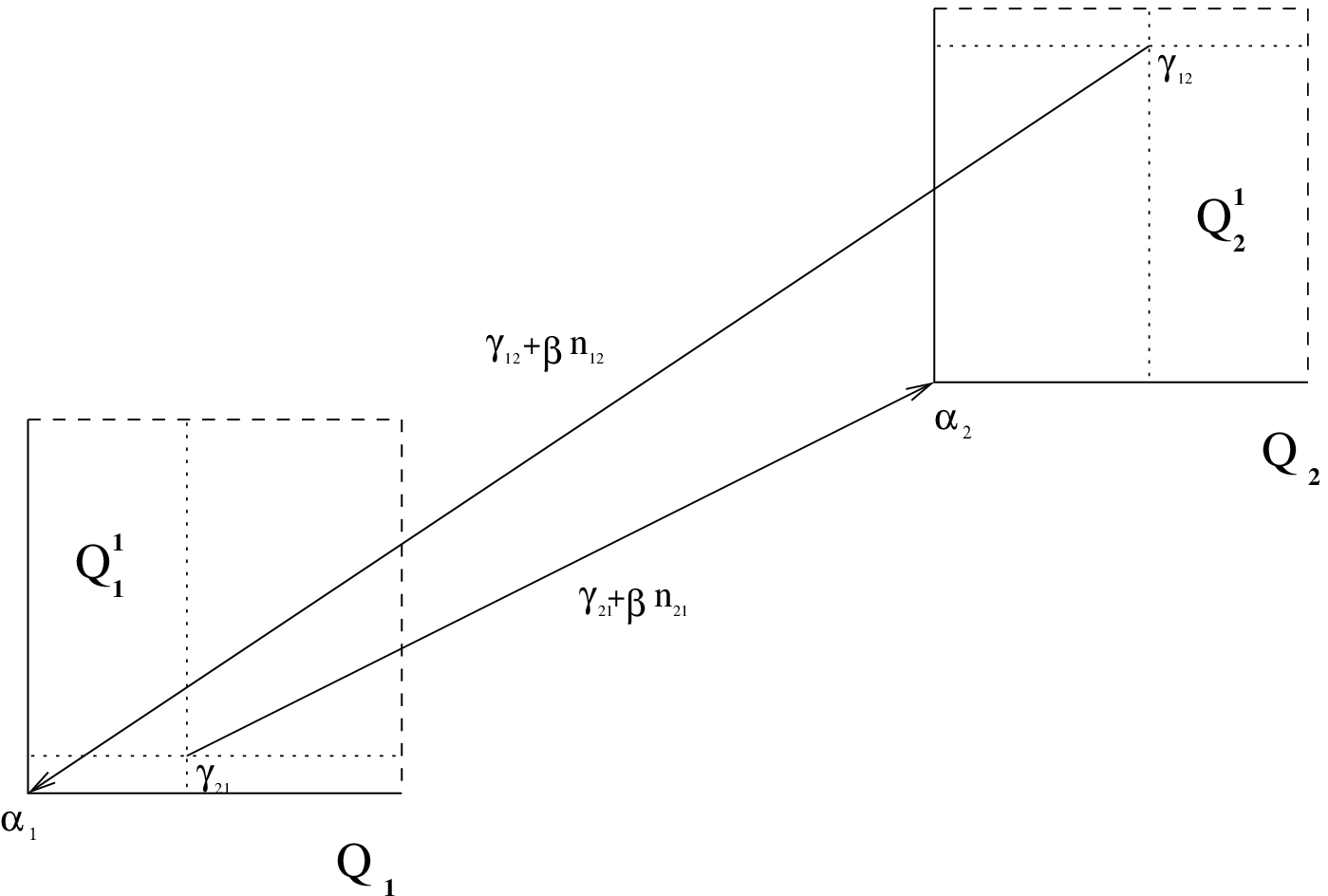}
 \end{center}
 \caption{ \label{dibuix}}
 \end{figure}

Although $\gamma_{jk}+\beta n_{jk}=\a_{k}\in Q_{k}$ and $\gamma_{jk}\in Q_{j},$
it can happen that for some $Q_{j}^{s}\subset Q_{j}$
$$(Q_{j}^{s}+\beta n_{jk})\cap Q_{k}=\emptyset,$$
see figure \ref{dibuix}.
To solve this difficulty we choose a correction
$$C_{js}(k)\in \{ -1,0,1 \}^{d}$$
in such a way that
$$Q_{j}^{s}+\beta (n_{jk}+C_{js}(k))= Q_{k}^{l_{s}},$$
for some $l_{s}\in \{1,\dots ,p^{d} \}.$

This shows that the intervals that form the cubes are always the same, but put in different orders,
i.e. for any $Q_{1}^{s}\subset Q_{1}$ we have a unique
rectangle in $Q_{j}$ that is translated by $\beta \Z^{d}$ in $Q_{1}^{s}.$
We call this rectangle
$Q_{j}^{s}.$ Then the intervals on each $Q_{j}$ are named depending on their translate on $Q_{1}.$ Using this
convention we have for $s=1,\dots ,p^{d}$ $k=1,\dots ,p$
$$Q_{1}^{s}+\beta(n_{1k}+C_{1s}(k))=Q_{k}^{s},$$
and
\begin{equation}                                                                        \label{relacio}
n_{kj}+C_{ks}(j)=n_{1j}-n_{1k}+C_{1s}(j)-C_{1s}(k).
\end{equation}
This is so because
the translation by $\beta\Z^{d}$ that moves $Q_{j}^{s}$ to their counterpart in $Q_{k}$ is unique.

Let
$$C_{js}=\{n_{jk}+C_{js}(k)  \}_{k=1,\dots p}$$
and notice that
$0=n_{jj}+C_{js}(j)\in C_{js}$ and as the cubes are disjoint the elements in $C_{js}$ are all different.

Now we have, for $\omega \in Q_{j}^{s}$
$$\omega +\beta n\in E \;\;\; \mbox{iff}\;\;\; n\in C_{js},$$
then for $\omega \in Q_{j}^{s}$ equation (\ref{marca}) becomes
\begin{eqnarray*}
    f(\omega) & = & \left( \frac{\beta}{\sqrt{2\pi}}\right)^{d}\sum_{l=1}^{p}s_{l}(\omega)
\sum_{n\in C_{js}}f(\omega+\beta n)e^{i\beta n\cdot k_{l}} \\
& = & \left( \frac{\beta}{\sqrt{2\pi}}\right)^{d}\sum_{n\in C_{js}}f(\omega+\beta n)
\sum_{l=1}^{p}s_{l}(\omega)
e^{i\beta n\cdot k_{l}}.
\end{eqnarray*}
If we find for all $j$ and $s$ a solution $x_{l}=x_{l}(j,s)$ $l=1,\dots ,p$ of the
linear system
\begin{equation}                                                            \label{system}
\left\{ \begin{array}{ll}
    \sum_{l=1}^{p}x_{l}e^{i\beta m\cdot k_{l}}=0, & \mbox{if}\;\; m\in C_{js}\setminus \{0\}
    \\
    \sum_{l=1}^{p}x_{l}=\left( \frac{\beta}{\sqrt{2\pi}}\right)^{-d}, & \mbox{otherwise}
\end{array}\right.
\end{equation}
we have finished, because defining $s_{l}(\omega)=x_{l}(j,s)$ for $\omega \in Q_{j}^{s}$ and
$$s_{l}(\omega)=\sum_{j=1}^{p}\sum_{s=1}^{p^{d}}x_{l}(j,s)\chi_{Q_{j}^{s}}(\omega)$$
we have (\ref{marca}).

We claim that
there exist $k_{l}\in \R^{d}$ for $l=1,\dots ,p$ such that
(\ref{system}) has a solution for all
$j=1,\dots ,p$ and $s=1,\dots ,p^{d}.$

Indeed,
letting $z_{ls}=e^{i\beta k_{l}^{s}}$ for $l=1,\dots ,p,$ $s=1,\dots ,d$
and $z_{l}=(z_{l1},\dots ,z_{ld})$ the linear system in (\ref{system}) has solution
if and only if there exists $z_{l}\in \mathbb{T}^{d},$ $l=1,\dots
,p$ such that
\begin{equation}                                                        \label{anulla}
\mathcal{D}(z_{11},\dots ,z_{pd})=\prod_{j=1}^{p}\prod_{s=1}^{p^{d}}
\det (\{ z_{l}^{n_{jk}+C_{js}(k)} \}_{k,l=1,\dots ,p })\neq 0.
\end{equation}

Indeed, letting
$$m_{js}=(\min_{k=1,\dots ,p} n_{jk}^{1}+C_{js}^{1}(k), \dots ,\min_{k=1,\dots ,p} n_{jk}^{d}+C_{js}^{d}(k)),$$
then
$$\det (\{ z_{l}^{n_{jk}+C_{js}(k)} \}_{k,l=1,\dots ,p })=
\frac{pol_{js}(z_{11},\dots ,z_{pd})}{z_{1}^{m_{js}}\cdots z_{p}^{m_{js}}}$$
where $pol_{js}$ is a non-zero polynomial in $\C^{pd}.$ Now
the first column expansion of the determinant gives us
$$\sum_{k=1}^{p}(-1)^{k}z_{1}^{n_{jk}+C_{js}(k)}M_{1k}(z_{2},\dots ,z_{p}).$$
Since the elements in $C_{js}$ are different, an induction reasoning provides the conclusion.

Moreover, the zero set of $\mathcal{D}$ has zero Lebesgue measure in the polydisc $\mathbb{T}^{pd},$
because $\mathcal{D}$ is an entire function
\cite[p. 218]{R}.
Then for almost all $(k_{1},\dots ,k_{p}) \in [0,\frac{2\pi}{\beta})^{pd},$
and all $j$ and $s$
the linear system in (\ref{system})
has a solution.

So far we have shown that for
$F\in PW_{E}$ such that $\mathcal{F}F\in
\mathcal{C}_{c}^{\infty}(E),$ there exist
$k_{l}\in \R^{d}$ with
\begin{equation}                                                                    \label{cond}
\mathcal{D}(e^{i\beta  k_{1}^{1}},\dots ,e^{i\beta k_{p}^{d}})\neq
0,
\end{equation}
so that taking
\begin{equation}                                                                           \label{func}
S_{l}(t)=\mathcal{F}^{-1}s_{l}(t)=\sum_{j=1}^{p}\sum_{s=1}^{p^{d}}x_{l}(j,s)
\mathcal{F}^{-1}\chi_{Q_{j}^{s}}(t), \;\; l=1,\dots ,p,
\end{equation}
we have
\begin{equation}                                                                        \label{eq-fin}
F(t)=\sum_{l=1}^{p}\sum_{n\in
    \Z^{d}}F(Wn+k_{l})S_{l}(t-Wn-k_{l}).
\end{equation}

The following results allows us to extend this formula to all $F\in
PW_{E}.$


\begin{lemma}                                                                   \label{lem:3}
    Take $k_{l}\in \R^{d}$ such that (\ref{cond}) holds and
    functions $S_{l}$ as in (\ref{func}). Then, if $\a=\cup_{l=1}^{p} \a^{l},$ $\a^{l}\in \ell^{2},$
    the series
\begin{equation}                                                                \label{f-al}
    F[\a](t)=\sum_{l=1}^{p}\sum_{n\in \Z^{d}}\a_{n}^{l} S_{l}
    (t-W n-k_{l}),
\end{equation}
    converge in $L^{2}(\R^{d}).$
    Moreover $F[\a]$ is an injective operator: there exists $B>0$
    independent of $\a$ such that
\begin{equation}\label{des}
    \| F[\a] \|_{L^{2}(\R^{d})}\le B \| \a
    \|_{\ell^{2}}.
\end{equation}
\end{lemma}

\begin{corollary}                                                               \label{corol}
    Taking functions $S_{l}$ as in (\ref{func}),
    and $k_{l} \in \R^{d}$ such that (\ref{cond}) holds, then
\begin{equation*}
    S_{l}(Wn+k_{s}-k_{l})=\delta_{n,0}\delta_{l,s},\;\; l,s=1,\dots
    ,p.
\end{equation*}
\end{corollary}

The proofs are given in the next section. We can now prove our main result.


\begin{theorem}                                                                 \label{teo:2}
    Let $E\subset \R^{d}$ be a finite union of cubes
    as in (\ref{conj}). If $W=\frac{2\pi}{\beta}Id$
    and
    $k_{l}\in \R^{p}$ are such that (\ref{cond}) holds then
\begin{equation*}
    \Lambda(k_{1},\dots ,k_{p})=\cup_{l=1}^{p}
    \left\{
    Wn+k_{l}
    \right\}_{n \in\Z^{d}}
\end{equation*}
    is a complete interpolating sequence for $PW_{E},$ i.e. there
    exists constants $A,B>0$ such that
\begin{equation}                                                                     \label{sampling2}
    A \| F \|^{2} \le \sum_{l=1}^{p}\sum_{n\in
    \Z^{d}}|F(Wn+k_{l})|^{2}\le B \| F \|^{2},\;\;\; F\in PW_{E}.
\end{equation}
    Given any $\a=\cup_{l=1}^{p}\a^{l}$ with
    $\a^{l}\in \ell^{2}$ there exists a unique $F\in PW_{E}$ such
    that
\begin{equation}                                                                        \label{interp}
    F(Wn+k_{l})=\a_{n}^{l}, \;\; n\in \Z^{d},l=1,\dots ,p.
\end{equation}
    Moreover, taking functions $S_{l}$ as in
    (\ref{func}),
    the reconstruction formula
    (\ref{eq-fin}) holds for
    all $F\in PW_{E},$ with convergence in $L^{2}-$norm and
    uniform convergence on products of horizontal strips.
\end{theorem}

\proof
    By the preceding construction we know that (\ref{eq-fin}) holds for
    all $F\in PW_{E}$ such that $\mathcal{F}F\in \mathcal{C}^{\infty}_{c}(E).$
    Now by Lemma \ref{lem:3} the operator on $PW_{E}$ defined by
    $$F\longmapsto \sum_{l=1}^{p}\sum_{n\in
    \Z^{d}}F(Wn+k_{l})S_{l}(t-Wn-k_{l}),$$
    is bounded
    and the sampling formula (\ref{eq-fin}) holds for $PW_{E}.$

    Since $\Lambda(k_{1},\dots ,k_{p})$ is uniformly discrete we have the right inequality
    in (\ref{sampling2}).
    In order to prove the left inequality, notice that by
    Lemma \ref{lem:3} we have that
    $\sum \a_{n}^{l}S_{l}(\cdot -Wn-k_{l})$ converges for all $\a^{l}\in \ell^{2},$
    therefore
    there exist a constant $B>0$ such that for $F\in PW_{E}$
    $$\sum_{l,n}|\ll F, S_{l}(\cdot -Wn-k_{l}) \rr|^{2}\le B\| F \|^{2},$$
    see \cite[Lemma 3.2.1.]{O}.
    Then for $F\in PW_{E}$ we have
\begin{eqnarray*}
    \| F \|^{4} & \le & \left(\sum_{l,n}|F(Wn+k_{l})|^{2}\right)
    \left( \sum_{l,n}|\ll F ,S_{l}(\cdot -W n-k_{l}) \rr|^{2}\right)
    \\
    & \le & B \| F \|^{2}\left(\sum_{l,n}|F(Wn+k_{l})|^{2}\right).
\end{eqnarray*}
    Defining for $\a=\cup_{l=1}^{p} \a^{l},$ $\a^{l}\in \ell^{2},$
    $F[\a]\in PW_{E}$
    as in (\ref{f-al}),
    and using Corollary \ref{corol} we get the interpolation property (\ref{interp})
\begin{eqnarray*}
    F[\a](Wm+k_{s}) & = &
    \sum_{l=1}^{p} \sum_{n\in \Z^{d}}\a_{n}^{l}
    S_{l}(W(m-n)+k_{s}-k_{l})
    \\
    & = &
    \sum_{l=1}^{p}\sum_{n\in \Z^{d}}\a^{l}_{n}
    \delta_{m,n}
    \delta_{l,s}
    =
    \a_{m}^{s}.
\end{eqnarray*}
    Finally, the uniform convergence on products of horizontal strips
    $$\{z\in \C^{d}:|\Im z^{j}|\le C_{j},j=1,\dots ,p  \},$$
    follows from the Phragm\'en-Lindel$\ddot{\rm o}$f inequality, see \cite[p 84]{Y}.
\qed


\begin{corollary}                                                           \label{coro}
    Let $E\subset \R^{d}$ be a bounded set with no empty interior.
    For all $\e>0$ there exist uniformly separated sequences $\L^{\e},\L_{\e}\subset \R^{d}$
    such that $\L_{\e}$ is sampling for $PW_{E}$, $\L^{\e}$ is
    interpolating for $PW_{E}$ and
    $$0\le (2\pi)^{d}D(\L_{\e})-|E|<\e ,\;\;\; 0\le |E|-(2\pi)^{d}D(\L^{\e})<\e.$$
\end{corollary}

\proof
    It is enough to choose $E_{\e},$ $E^{\e}$
    as
    in (\ref{conj}) such that
    $E^{\e}\subset E \subset E_{\e}$ and $|E_{\e}|-|E|,|E|-|E^{\e}|<\e.$
    Now we take a complete interpolating sequence $\L^{\e}$ for $PW_{E^{\e}}$ (that obviously is also interpolating for
    $PW_{E}$) and $\L_{\e}$ a
    complete interpolating sequence for $PW_{E_{\e}}$ (that is also sampling for
    $PW_{E}$) and apply Landau's Theorem \ref{lan}.
\qed

\section{Proofs of technical results}                    \label{proofs}

In this section we will proof the technical
results used in proving
Theorem \ref{teo:1}.

\vskip .2in


\proof [Lemma \ref{lem:3}]
    Consider, given $M\in \N$
$$F_{M}[\a](t)=\sum_{l=1}^{p}\sum_{|n| \le M}\a_{n}^{l} S_{l}
    (t-W n-k_{l}).$$
    For all $f\in L^{2}(\R^{d})$
    we have
\begin{eqnarray*}
    \ll F_{M}[\a], f\rr & = & \sum_{l,j=1}^{p}\sum_{s=1}^{p^{d}}x_{l}(j,s)
    \sum_{|n| \le M}\a_{n}^{l} \ll e^{-i\omega\cdot (\frac{2\pi}{\beta}n+k_{l})},\chi_{Q_{j}^{s}}\widehat{f}\rr
    \\
    & = & \sum_{l,j=1}^{p}\sum_{s=1}^{p^{d}}x_{l}(j,s)
    \sum_{|n| \le M}\a_{n}^{l} \mathcal{P}_{PW_{Q_{j}^{s}}}(f)\left(\frac{2\pi}{\beta}n+k_{l}\right),
\end{eqnarray*}
    where $\mathcal{P}_{PW_{Q_{j}^{s}}}(f)=\mathcal{F}^{-1}(\chi_{Q_{j}^{s}}\mathcal{F}f)$ is the
    orthogonal projection of $f\in L^{2}(\R^{d})$ into $PW_{Q_{j}^{s}}.$
    Then, using Lemma \ref{lem:1} $(d),$ we obtain
\begin{eqnarray*}
    |\ll F_{M}[\a], f\rr | & \le &
    \sum_{l,j=1}^{p}\sum_{s=1}^{p^{d}}|x_{l}(j,s)|
    \| \a^{l} \|_{\ell^{2}} \| \mathcal{P}_{PW_{Q_{j}^{s}}}(f) \|_{L^{2}} \\
    & \le &
     B \| f \|_{L^{2}} \sum_{l=1}^{p}\| \a^{l} \|_{\ell^{2}} \sum_{j=1}^{p} \sum_{s=1}^{p^{d}}|x_{l}(j,s)|
     \| \mathcal{P}_{PW_{Q_{j}^{s}}}\| \\
    & \le &
    B \| f \|_{L^{2}}\sum_{l=1}^{p}\| \a^{l} \|_{\ell^{2}},
\end{eqnarray*}
    where the constant $B$ does not depend on $\a.$

Now we want to prove the injectivity part.
By taking Fourier transform in (\ref{f-al}) we obtain

\begin{equation}                                                                            \label{eq:2}
    0=\sum_{j=1}^{p}\sum_{s=1}^{p^{d}}\chi_{Q_{j}^{s}}(\omega)\sum_{l=1}^{p}x_{l}(j,s)
    e^{-i\omega \cdot k_{l}}
    \left\{ \sum_{n\in \Z^{d}}\a_{n}^{l}e^{-i\frac{2\pi}{\beta}\omega \cdot n}\right\}.
\end{equation}
    We call
    $$f_{l}(\omega)=\sum_{n\in \Z^{d}}\a_{n}^{l}
    e^{-i \frac{2\pi}{\beta}\omega\cdot n},\;\;\;l=1,\dots ,p.$$
    If we show that
$$
\begin{array}{ll}
f_{l}(\omega)=0, & l=1,\dots p,
\end{array}
$$
    we are done.
    The functions $f_{l}$ are $\beta-$periodic on each variable,
    hence it is enough
    to proof $f_{l}|Q_{1}=0.$
We will prove this on each $Q_{1}^{s}$ separately. The key is
property (\ref{relacio}).

Let us take $\omega \in Q_{1}^{s}.$ Then, as $f_{l}(\omega+\beta(n_{1k}+C_{1s}(k)))=f_{l}(\omega)$
    and $\omega+\beta(n_{1k}+C_{1s}(k))\in Q_{ks},$ substituting into (\ref{eq:2}) we obtain, for $k=1,\dots ,p$

\begin{eqnarray*}
0 & = & \sum_{l=1}^{p}x_{l}(j,s)
    e^{-i(\omega+\beta(n_{1k}+C_{1s}(k))) \cdot k_{l}}f_{l}(\omega)
    \\
    & = &
    \sum_{l=1}^{p}x_{l}(j,s)
    e^{-i\beta(n_{1k}+C_{1s}(k)) \cdot k_{l}}[e^{-i\omega \cdot k_{l}}f_{l}(\omega)].
\end{eqnarray*}
    We have a linear system with unknowns $e^{-i\omega \cdot k_{l}}f_{l}(\omega)$
    and
    coefficient's matrix
$$(x_{l}(k,s) e^{-i\beta (n_{1k}+C_{1s}(k)) \cdot k_{l}})_{k,l}=(a_{kl})_{k,l}.$$

We will see that this system has only the trivial solution by showing that it has invertible coefficient's matrix.
This point becomes quite easy with the notation we have adopted.
Indeed, the linear system (\ref{system}) that we have used to determine $s_{l}$ in $Q_{1}^{s}$ is given by
$$(e^{i\beta (n_{1j}+C_{1s}(j)) \cdot k_{l}})_{l,j}=(b_{lj})_{l,j}.$$
Now by (\ref{relacio}), the elements of the matrix product are
\begin{eqnarray*}
\sum_{l=1}^{p}a_{kl}b_{lj} & = & \sum_{l=1}^{p}x_{l}(k,s)
e^{i\beta (n_{1j}-n_{1k}+C_{1s}(j)-C_{1s}(k)) \cdot k_{l}}
\\
& = &
\sum_{l=1}^{p}x_{l}(k,s)
e^{i\beta (n_{kj}+C_{ks}(j)) \cdot k_{l}}=\delta_{kj}\left( \frac{\sqrt{2\pi}}{\beta}\right)^{d}
\end{eqnarray*}
which is the desired conclusion.
\qed


\proof [Corollary \ref{corol}] Let $s\in\{ 1,\dots ,p\}$ be fixed.
The sampling result applied to $S_{s}(\cdot-k_{l})\in PW_{E}$ yields
\begin{equation*}
    S_{s}(t-k) =\sum_{l=1}^{p}\sum_{n\in \Z^{d}}
    S_{s}(Wn+k_{l}-k_{s})
    S_{l}(t-Wn-k_{l}).
\end{equation*}
    Now applying Lemma \ref{lem:3} to
\begin{equation*}
    \sum_{l=1}^{p}\sum_{n\in \Z^{d}}  \left(
    \delta_{n,0}\delta_{l,s}-
    S_{s}(Wn+k_{l})
    \right)
    S_{l}(t-Wn-k_{l})=0,
\end{equation*}
    we obtain the result.
\qed
\vskip .2in


\end{document}